\documentclass[12pt,british]{article}
\usepackage[T1]{fontenc}
\usepackage{geometry}
\usepackage{color}
\usepackage{amsmath}
\usepackage{amsthm}
\usepackage{amssymb}
\PassOptionsToPackage{normalem}{ulem}
\usepackage{ulem}

\makeatletter
\theoremstyle{plain}
\newtheorem{thm}{\protect\theoremname}
\theoremstyle{plain}
\newtheorem{prop}[thm]{\protect\propositionname}

\usepackage{textcomp}

\usepackage{babel}
\addto\captionsbritish{}
  \addto\captionsbritish{}
  \addto\captionsbritish{\renewcommand{\propositionname}{Proposition}}
  \addto\captionsbritish{}
  \addto\captionsbritish{\renewcommand{\theoremname}{Theorem}}
  \addto\captionsenglish{}
  \addto\captionsenglish{}
  \addto\captionsenglish{\renewcommand{\propositionname}{Proposition}}
  \addto\captionsenglish{}
  \addto\captionsenglish{\renewcommand{\theoremname}{Theorem}}

  \providecommand{\propositionname}{Proposition}
  
\providecommand{\theoremname}{Theorem}

\usepackage{babel}
\addto\captionsbritish{}
\addto\captionsbritish{\renewcommand{\propositionname}{Proposition}}
\addto\captionsbritish{\renewcommand{\theoremname}{Theorem}}
\addto\captionsenglish{}
\addto\captionsenglish{\renewcommand{\propositionname}{Proposition}}
\addto\captionsenglish{\renewcommand{\theoremname}{Theorem}}

\providecommand{\propositionname}{Proposition}
\providecommand{\theoremname}{Theorem}

\makeatother

\usepackage{babel}
\providecommand{\propositionname}{Proposition}
\providecommand{\theoremname}{Theorem}

\begin{document}
\title{Ricci solitons with an orthogonally intransitive $2$-dimensional
Abelian Killing algebra}
\author{Diego Catalano Ferraioli$^{1}$}

\maketitle
\footnotetext[1]{Instituto de Matem\'{a}tica e Estat\'{\i}stica-
Universidade Federal da Bahia, Campus de Ondina, Av. Adhemar de Barros,
S/N, Ondina - CEP 40.170.110 - Salvador, BA - Brazil, e-mail: diego.catalano@ufba.br.}
\begin{abstract}
In this paper we report on a local classification of four dimensional
Ricci solitons which have a $2$-dimensional Abelian Killing algebra
$\mathcal{G}_{2}$, whose Killing leaves are non-null and orthogonally
intransitive. The classification is obtained under the following additional
assumptions: \textit{(i)} the curvature vector field, of the submersion
defined by $\mathcal{G}_{2}$, is a null vector field; \textit{(ii)
}$\mathcal{G}_{2}$ has a null vector; \textit{(iii) }the vector field
of the Ricci soliton is tangent to the Killing leaves and a symmetry
of the orthogonal distribution. Since there are only few examples
of orthogonally intransitive Einstein metrics, and even less is known
about orthogonally intransitive Ricci solitons, we believe that these
results can help fill this gap in the literature.

\vspace{0.3cm}

\noindent Mathematics Subject Classification 53B20, 53C21, 53C50, 53C25  
\end{abstract}

\section{Introduction}

Ricci solitons were introduced by Richard Hamilton in \cite{Hamilton},
as self-similar solutions of the Ricci flow. Since then, they have
been a subject of intense research, especially after they appeared
to be fundamental in the proof of the Poincar\'{e} conjecture (see
\cite{CLN,CetAl,Hamilton} and references therein).

We recall that a pseudo-Riemannian manifold $(\mathcal{M},\mathbf{g})$
is a Ricci soliton if there exists a vector field $X$ on $\mathcal{M}$
and a real constant $\Lambda$ such that 
\begin{equation}
\text{Ric}(\mathbf{g})+\frac{1}{2}L_{X}g=\Lambda g,\label{eq:Ricci-soliton}
\end{equation}
where $Ric(\mathbf{g})$ is the Ricci tensor of $\mathbf{g}$ and
$L_{X}(\mathbf{g})$ is the Lie derivative of $\mathbf{g}$ with respect
to $X$. In particular, the Ricci soliton is said to be expanding,
steady or shrinking if $\Lambda<0$, $\Lambda=0$ or $\Lambda>0$,
respectively. 

As such the notion of Ricci soliton is a generalization of that of
Einstein manifold. Therefore, it is natural to think that the study
and classification of Ricci solitons, like that of Einstein metrics,
may significantly depend on the properties of their Killing algebra
$\mathfrak{Kill}(\mathbf{g})$.

In four dimensional case, especially in physics \cite{S-K-M-H-H},
a rather rich class of examples of exact solutions of the equation
$\text{Ric}(\mathbf{g})=\Lambda g$ is that with a $2$-dimensional
Abelian Killing algebra $\mathcal{G}_{2}:=<\xi_{(1)},\xi_{(2)}:\,[\xi_{(1)},\xi_{(2)}]=0>$.
In that class, however, the vast majority of known solutions possess
the further property of being orthogonally transitive, that is the
distribution $\Xi:=\mathit{span}\{\xi_{(1)},\xi_{(2)}\}$ admits a
completely integrable orthogonal distribution $\Xi^{\perp}$. In fact,
only a handful of orthogonally intransitive metrics are known.  In
the paper \cite{Cat-Mar} we found new orthogonally intransitive metrics
under the additional assumption that some first order differential
invariants are identically zero.

Based on \cite{Cat-Mar} we recently studied also the case of four
dimensional Ricci solitons which have a $2$-dimensional Abelian Killing
algebra $\mathcal{G}_{2}$, whose Killing leaves are non-null and
orthogonally intransitive. As a result we obtained a local classification
under the following additional assumptions:
\begin{description}
\item [{\textit{(i)}}] the curvature vector field $\mathcal{C}$, of the
submersion defined by $\mathcal{G}_{2}$, is a null vector field; 
\item [{\textit{(ii)}}] $\mathcal{G}_{2}$ has a null vector;\textsl{ }
\item [{\textit{(iii)}}] $X$ is tangent to the Killing leaves, i.e., $X\in\Xi$,
and is a symmetry of $\Xi^{\perp}$.
\end{description}
In this short note we will summarise the main results of this classification,
leaving out the detailed proofs, which will be left to another paper. 

\section{Preliminaries\label{sec:Preliminaries}}

Given a pseudo-Riemannian\textcolor{black}{{} manifold} $(\mathcal{M},\mathbf{g})$,
whose Killing algebra is a $2$-dimensional Abelian algebra of vector
fields $\mathcal{G}_{2}:=<\xi_{(1)},\xi_{(2)}:\,[\xi_{(1)},\xi_{(2)}]=0>$,
one can always choose local coordinates $\{t_{1},t_{2},z_{1},z_{2}\}$
on $\mathcal{M}$ such that $\xi_{(i)}=\partial_{z_{i}}$ and the
integral manifolds of the completely integrable distribution $\Xi:=\mathit{span}\{\xi_{(1)},\xi_{(2)}\}$
are surfaces characterised by the constancy of $t_{1}$ and $t_{2}$.
We will call Killing leaves the integral manifolds of $\Xi$, and
will refer to such a kind of coordinates $\{t_{1},t_{2},z_{1},z_{2}\}$
as local adapted coordinates. 

Adapted coordinates are obviously not unique, and one can consider
the Lie pseudogroup $\mathfrak{G}$ of adapted coordinates transformations.
By definition, these are coordinate transformations $\bar{t}_{i}=\bar{t}_{i}(t_{1},t_{2},z_{1},z_{2})$,
$\bar{z}_{i}=\bar{z}_{i}(t_{1},t_{2},z_{1},z_{2})$ such that $\mathcal{G}_{2}$
is generated by $\partial/\partial\bar{z}_{i}$, $i=1,2$, and the
leaves of $\Xi$ are surfaces characterised by the constancy of $\bar{t}_{1}$
and $\bar{t}_{2}$. It can be shown (see \cite{Cat-Mar}) that these
transformations have the form 
\begin{equation}
\bar{t}_{i}=\phi_{i}(t_{1},t_{2}),\qquad\bar{z}_{i}=\alpha_{i}^{j}\,z_{j}+\psi_{i}(t_{1},t_{2}),\label{eq:adapt_coord_transf}
\end{equation}
where $(\alpha_{i}^{j})\in\mathrm{GL}(2,\mathbb{R})$ and $\phi_{i}(t_{1},t_{2})$,
$\psi_{i}(t_{1},t_{2})$ are arbitrary differentiable functions such
that 
\[
J_{\phi}=\left|\begin{array}{@{}ll@{}}
\partial_{t_{1}}\phi_{1} & \partial_{t_{2}}\phi_{1}\\
\partial_{t_{1}}\phi_{2} & \partial_{t_{2}}\phi_{2}
\end{array}\right|\neq0.
\]
The Einstein summation convention of summing over repeated indices
is used here and throughout the paper. Thus, in adapted coordinates
the metric $\mathbf{g}$ takes the form 
\[
\mathbf{g}=b_{ij}(t_{1},t_{2})\,dt_{i}\,dt_{j}+2f_{ik}(t_{1},t_{2})\,dt_{i}\,dz_{k}+h_{kl}(t_{1},t_{2})\,dz_{k}dz_{l},
\]
where $b_{21}=b_{12}$ and $h_{21}=h_{12}$. In particular $h_{kl}(t_{1},t_{2})\,dz_{k}dz_{l}$
describes the restriction of $\mathbb{\mathbf{g}}$ to the Killing
leaves. Hence the Killing leaves are not null iff $\det(h_{ij})\neq0$
everywhere. In \cite{Cat-Mar}, however, we found it very convenient
to rewrite $\mathbf{g}$ in the following alternative form 
\begin{equation}
\mathbf{g}=\tilde{g}_{ij}(t_{1},t_{2})\,dt_{i}\,dt_{j}+h_{kl}(t_{1},t_{2})\,\left(dz_{k}+f_{i}^{k}dt_{i}\right)\left(dz_{l}+f_{j}^{l}dt_{j}\right)\label{eq:Geroch}
\end{equation}
where 
\[
\tilde{g}_{ij}:=b_{ij}-f_{ik}f_{jl}h^{kl},\qquad\qquad f_{j}^{k}:=f_{js}h^{sk},
\]
and $h^{kl}$ denote the elements of the inverse matrix $(h_{ij})^{-1}$. 

An important advantage of (\ref{eq:Geroch}) is that $\tilde{\mathbf{g}}=\tilde{g}_{ij}(t_{1},t_{2})\,dt_{i}\,dt_{j}$
defines locally a natural metric on the orbit space $\mathcal{S}=\mathcal{M}/G_{2}$,
where by $G_{2}$ we mean the Lie group of transformations generated
by $\mathcal{G}_{2}$ on $\mathcal{M}$. Moreover, it can be checked
that the projection $\pi:\mathcal{M}\rightarrow\mathcal{S}$, $(t_{1},t_{2},z_{1},z_{2})\mapsto(t_{1},t_{2})$,
is locally a pseudo-Riemannian submersion, with respect to the metrics
$\mathbf{g}$ and $\tilde{\mathbf{g}}$. Throughout the paper we will
refer to $\tilde{\mathbf{g}}$ as the orbit metric.

Then, in view of (\ref{eq:Geroch}), one can easily check that $\Xi^{\perp}$
is generated by the vector fields 
\[
\mathbf{e}_{j}=\frac{\partial}{\partial t_{j}}-f_{j}^{1}\frac{\partial}{\partial z_{1}}-f_{j}^{2}\frac{\partial}{\partial z_{2}},\qquad j=1,2.
\]
Moreover, the curvature vector field $\mathcal{C}$ of the submersion
can be written as
\[
\mathcal{C}={\displaystyle \frac{\partial_{t_{2}}f_{1}^{1}-\partial_{t_{1}}f_{2}^{1}}{\sqrt{\left|\det\tilde{\mathbf{g}}\right|}}}\frac{\partial}{\partial z_{1}}+{\displaystyle \frac{\partial_{t_{2}}f_{1}^{2}-\partial_{t_{1}}f_{2}^{2}}{\sqrt{\left|\det\tilde{\mathbf{g}}\right|}}}\frac{\partial}{\partial z_{2}}
\]
 and its squared length reads
\[
\mathbf{g}(\mathcal{C},\mathcal{C})=h_{kl}\mathcal{C}^{k}\mathcal{C}^{l}={\displaystyle \frac{h_{kl}(\partial_{t_{2}}f_{1}^{k}-\partial_{t_{1}}f_{2}^{k})(\partial_{t_{2}}f_{1}^{l}-\partial_{t_{1}}f_{2}^{l})}{\left|\det\tilde{\mathbf{g}}\right|}}.
\]

Throughout the paper we will say that the metric or even the Killing
leaves are orthogonally intransitive if $\Xi^{\perp}$ is not integrable.
The reader can easily check that $\Xi^{\perp}$ is integrable iff
$\mathcal{C}=0$. Hence, the orthogonally intransitive case is characterised
by the condition $\mathcal{C}\neq0$.

\section{Classification results\label{sec:Classification-results}}

We assume that $\mathbf{g}$ is a pseudo-Riemannian metric with a
$2$-dimensional Abelian Killing algebra $\mathcal{G}_{2}$, whose
Killing leaves are non-null and orthogonally intransitive. In addition
to this, we also assume the following:
\begin{description}
\item [{\textit{(i)}}] $\mathbf{g}(\mathcal{C},\mathcal{C})=0$; 
\item [{\textit{(ii)}}] $\mathcal{G}_{2}$ has a null vector;\textsl{ }
\item [{\textit{(iii)}}] $X\in\Xi$ and is a symmetry of $\Xi^{\perp}$.
\end{description}
In view of assumptions \textit{(i)}-\textit{(iii)}, one can always
choose the adapted coordinates $\{t_{1},t_{2},z_{1},z_{2}\}$ in such
a way that $\xi_{(2)}=\partial_{z_{2}}$ is light-like and in addition
\[
\left(\tilde{g}_{ij}\right)=\left(\begin{array}{cc}
e^{-2P}\qquad & 0\vspace{5pt}\\
0\qquad & \pm e^{-2P}
\end{array}\right),\qquad\left(h_{ij}\right)=\left(\begin{array}{cc}
h_{11}\qquad & h_{12}\vspace{5pt}\\
h_{12}\qquad & 0
\end{array}\right).
\]
Moreover, in view of (\ref{eq:adapt_coord_transf}), it is also possible
to further adapt the coordinates in such a way that $f_{2}^{1}=f_{2}^{2}=0$.

Thus, without loss of generality, one can assume that 
\begin{equation}
\mathbf{g}=e^{-2P}\,\left(dt_{1}^{2}+\epsilon_{0}dt_{2}^{2}\right)+h_{11}\,\left(dz_{1}+f_{1}^{1}dt_{1}\right)^{2}+2h_{12}\,\left(dz_{1}+f_{1}^{1}dt_{1}\right)\left(dz_{2}+f_{1}^{2}dt_{1}\right),\label{eq:Forma_soluz}
\end{equation}
with $\epsilon_{0}=\pm1$ and $P=P(t_{1},t_{2})$, $h_{ij}=h_{ij}(t_{1},t_{2})$,
$f_{j}^{j}=f_{i}^{j}(t_{1},t_{2})$ differentiable functions satisfying
the non degeneracy condition $h_{12}\neq0$ and also the condition
$\mathcal{C}\neq0$, i.e.,
\[
\left(f_{1,t_{2}}^{1}\right)^{2}+\left(f_{1,t_{2}}^{2}\right)^{2}\neq0.
\]
We found that, under above assumptions, the analysis of equation (\ref{eq:Ricci-soliton})
naturally splits into the following two main cases:
\begin{description}
\item [{(I)}] $\underline{f_{1,t_{2}}^{1}\neq0}$;
\item [{(II)}] $\underline{f_{1,t_{2}}^{1}=0}$.
\end{description}
In particular, the case (II) further splits into the three sub-cases:
\begin{description}
\item [{(II-1)}] $\underline{h_{12,t_{1}}\neq0,\;\Lambda\neq0}$; 
\item [{(II-2)}] $\underline{h_{12,t_{1}}\neq0,\;\Lambda=0}$;
\item [{(II-3)}] $\underline{h_{12,t_{1}}=0}$.
\end{description}
The complete and detailed proofs of the results of our study will
be published in another paper, here we limit ourselves to give an
account of the obtained results. 

\subsection*{Case I}

Here is the main result when $f_{1,t_{2}}^{1}\neq0$.
\begin{thm}
For a suitable choice of adapted coordinates $\{t_{1},t_{2},z_{1},z_{2}\}$,
any Ricci soliton with a $2$-dimensional Abelian Killing algebra
$\mathcal{G}_{2}$, which is orthogonally intransitive, has non-null
Killing leaves and satisfies conditions (i), (ii) and (iii) with $f_{1,t_{2}}^{1}\neq0$,
can be written in the form 
\[
g=e^{-2P}\left(dt_{1}^{2}+\epsilon_{0}dt^{2}\right)-2\left(\int e^{-2P}dt_{2}\right)\left(c\,dt_{1}dz_{1}-dt_{1}dz_{2}\right)-2c\,dz_{1}^{2}+2dz_{1}dz_{2}
\]
with
\[
X=\left(\frac{\epsilon_{0}\,c}{2}z_{1}-\frac{\epsilon_{0}}{2}z_{2}+a_{1}\right)\partial_{z_{1}}+\left(\left(\frac{\epsilon_{0}\,c^{2}}{2}-2\Lambda\,c\right)z_{1}+\left(2\Lambda-\frac{\epsilon_{0}\,c}{2}\right)z_{2}+a_{2}\right)\partial_{z_{2}},
\]
where $\Lambda,c,a_{i}\in\mathbb{R}$ are arbitrary constants, $\epsilon_{0}=\pm1$
and in addition $P=P(t_{1},t_{2})$ is a differentiable function satisfying
the Liouville type differential equation 
\[
P_{,t_{1}t_{1}}+\epsilon_{0}P_{,t_{2}t_{2}}=\Lambda e^{-2P}.
\]
In particular, the Ricci soliton can be either expanding ($\Lambda<0$),
steady ($\Lambda=0$) or shrinking ($\Lambda>0$), moreover the orbit
metric $\tilde{\mathbf{g}}=e^{-2P}(dt_{1}^{2}+\epsilon_{0}dt_{2}^{2})$
has constant Gauss curvature $\Lambda$.
\end{thm}

\subsection*{Case II-1}

In this case one finds that $h_{12}=R$, with $R=R(t_{1})$ a positive
valued differentiable function such that $R'=R^{3/2}+c$, with $c\in\mathbb{R}$.
Thus, in view of this, one is naturally lead to distinguish the two
cases $c\neq0$ and $c=0$. This way one gets the following 
\begin{prop}
For a suitable choice of adapted coordinates $\{t_{1},t_{2},z_{1},z_{2}\}$,
any Ricci soliton with a $2$-dimensional Abelian Killing algebra
$\mathcal{G}_{2}$, which is orthogonally intransitive, has non-null
Killing leaves and satisfies conditions (i), (ii) and (iii) with $f_{1,t_{2}}^{1}=0$
and $h_{12,t_{1}}\neq0$, $\Lambda\neq0$, reduces to one of the following
two types:
\begin{description}
\item [{\uline{Type\mbox{$\,$}A'}$\,$$(c\neq0)$:}] 
\[
\begin{array}{l}
\mathbf{g}=-\frac{3\left(R^{3/2}+c\right)}{4\Lambda\sqrt{R}}\,\left(dt_{1}^{2}+\epsilon_{0}dt_{2}^{2}\right)-\frac{3\,t_{2}}{\Lambda}\left(1+\frac{c}{R^{3/2}}\right)\,dt_{1}dz_{1}+\frac{4}{3c\Lambda}\left(\frac{R^{2}S-\epsilon_{0}\sqrt{R}}{R}\right)dz_{1}^{2}+2Rdz_{1}dz_{2}\vspace{10pt}\end{array}
\]
with
\[
X=a_{1}\partial_{z_{1}}+\left(\frac{4a_{2}}{3c\Lambda}z_{1}+a_{3}\right)\partial_{z_{2}},
\]
where $\epsilon_{0}=\pm1$ and $\Lambda,c,a_{i}\in\mathbb{R}$ are
arbitrary constants satisfying $\Lambda<0$, $c\neq0$, moreover $R=R(t_{1})>0$
and $S=S(t_{1},t_{2})$ are differentiable functions satisfying the
differential equations 
\[
R'=R^{3/2}+c,\qquad\qquad S_{,t_{1}t_{1}}+\epsilon_{0}S_{,t_{2}t_{2}}+\frac{\left(R^{3/2}+c\right)}{R}S_{,t_{1}}=-\frac{3a_{2}\left(R^{3/2}+c\right)}{2\Lambda\sqrt{R}}
\]
and the inequality $R^{3/2}+c>0$. \\
In particular, the Ricci soliton is expanding and the orbit metric
$\tilde{\mathbf{g}}=-(3/4)\left(R^{3/2}+c\right)\Lambda^{-1}R^{-1/2}\,(dt_{1}^{2}+\epsilon_{0}dt_{2}^{2})$
has Gauss curvature 
\[
K=\frac{\Lambda\left(R^{5}+c^{3}\sqrt{R}+3c^{2}R^{2}+3cR^{7/2}\right)}{3\left(R^{5/2}+cR\right)^{2}}.
\]
\item [{\uline{Type\mbox{$\,$}B'}$\,$$(c=0)$:}] 
\[
\mathbf{g}=-\frac{3R}{4\Lambda}\,\left(dt_{1}^{2}+\epsilon_{0}dt_{2}^{2}\right)-\frac{3\,t_{2}}{\Lambda}\,dt_{1}dz_{1}+SRdz_{1}^{2}+2Rdz_{1}dz_{2}
\]
with
\[
X=a_{1}\partial_{z_{1}}+\left(a_{2}z_{1}+a_{3}\right)\partial_{z_{2}},
\]
where $\epsilon_{0}=\pm1$ and $\Lambda,a_{i}\in\mathbb{R}$ are arbitrary
constants satisfying $\Lambda<0$, moreover $R=R(t^{1})>0$ and $S=S(t^{1},t^{2})$
are differentiable functions satisfying the differential equations
\[
R'=R^{3/2},\qquad\qquad S_{,t_{1}t_{1}}+\epsilon_{0}S_{,t_{2}t_{2}}+\sqrt{R}\,S_{,t_{1}}=-\frac{3}{2\Lambda}\left(a_{2}R+\frac{2\epsilon_{0}}{R^{2}}\right).
\]
In particular, the Ricci soliton is expanding and the orbit metric
$\tilde{\mathbf{g}}=-(3/4)\Lambda^{-1}R\,(dt_{1}^{2}+\epsilon_{0}dt_{2}^{2})$
has constant Gauss curvature $\Lambda/3$.
\end{description}
\end{prop}

\bigskip{}
Thus, in view of this proposition, one can further refine the choice
of the adapted coordinates, and perform the following coordinate transformations:
\begin{itemize}
\item for Type A' 
\[
t_{1}\mapsto\epsilon\sqrt{R(t_{1})},\qquad t_{2}\mapsto-\epsilon\,\frac{t_{2}}{2c_{1}}\,\sqrt{-\frac{3}{\Lambda}},\qquad z_{1}\mapsto2\epsilon\,c_{1}\,z_{1}\,\sqrt{-\frac{1}{3\Lambda}},\qquad z_{2}\mapsto-\epsilon\,z_{2}\,\frac{\sqrt{-3\Lambda}}{2\,c_{1}},
\]
where $\epsilon=c/|c|$ and $c_{1}=1/\sqrt{\epsilon\,c}$;
\item for Type B'
\[
t_{1}\mapsto\frac{1}{\sqrt{R(t_{1})}}\sqrt[3]{\frac{3}{\Lambda}},\qquad t_{2}\mapsto\frac{t_{2}}{2}\sqrt[3]{\frac{3}{\Lambda}},\qquad z_{1}\mapsto2\,z_{1}\,\sqrt[3]{\frac{3}{\Lambda}},\qquad z_{2}\mapsto z_{2}\,\sqrt[3]{\frac{3}{\Lambda}}.
\]
\end{itemize}
This way, by further adapting the arbitrary function $S$ and the
remaining constants, one can prove the following
\begin{thm}
For a suitable choice of adapted coordinates $\{t_{1},t_{2},z_{1},z_{2}\}$,
any Ricci soliton with a $2$-dimensional Abelian Killing algebra
$\mathcal{G}_{2}$, which is orthogonally intransitive, has non-null
Killing leaves and satisfies conditions (i), (ii) and (iii) with $f_{1,t_{2}}^{1}=0$
and $h_{12,t_{1}}\neq0$, $\Lambda\neq0$, reduces to one of the following
two types:
\begin{description}
\item [{\uline{Type\mbox{$\,$}A}:}] 
\[
\begin{array}{l}
\mathbf{g}=-{\displaystyle \frac{3\,c_{1}^{2}t_{1}}{\Lambda\left(c_{1}^{2}t_{1}^{3}+1\right)}}\,dt_{1}^{2}+{\displaystyle \frac{\epsilon_{0}\left(c_{1}^{2}t_{1}^{3}+1\right)}{t_{1}}}\,dt_{2}^{2}-{\displaystyle \frac{6\,\epsilon\,t_{2}}{t_{1}^{2}}}\,dt_{1}dz_{1}+{\displaystyle \frac{\epsilon\,\psi\,t_{1}^{3}+\epsilon_{0}}{t_{1}}}\,dz_{1}^{2}-2t_{1}^{2}\,dz_{1}dz_{2}\vspace{10pt}\end{array}
\]
with
\[
X=a_{1}\partial_{z_{1}}+\left(a_{2}z_{1}+a_{3}\right)\partial_{z_{2}},
\]
where $\epsilon_{0},\epsilon=\pm1$ and $\Lambda,c_{1},a_{i}\in\mathbb{R}$
are arbitrary constants satisfying $\Lambda<0$, $c_{1}\neq0$, moreover
$\psi=\psi(t_{1},t_{2})$ is a differentiable function satisfying
the differential equation 
\[
\psi_{,t_{1}t_{1}}-\frac{3\,\epsilon_{0}\,\epsilon\,c_{1}\,t_{1}^{2}}{\Lambda\,(\epsilon\,t_{1}^{3}+c_{1})^{2}}\psi_{,t_{2}t_{2}}+\frac{\left(4\,\epsilon\,t_{1}^{3}+c_{1}\right)}{t_{1}\left(\epsilon\,t_{1}^{3}+c_{1}\right)}\psi_{,t_{1}}=\frac{6a_{2}t_{1}}{\Lambda\,\left(\epsilon\,t_{1}^{3}+c_{1}\right)}.
\]
In particular the Ricci soliton is expanding and the orbit metric
$\tilde{\mathbf{g}}=-3\,c_{1}^{2}t_{1}\Lambda^{-1}\left(c_{1}^{2}t_{1}^{3}+1\right)^{-1}dt_{1}^{2}+\epsilon_{0}\left(c_{1}^{2}t_{1}^{3}+1\right)t_{1}^{-1}dt_{2}^{2}$
has Gauss curvature 
\[
K=\frac{\epsilon\,\Lambda\,\left(\epsilon\,t_{1}^{3}+c_{1}\right)}{3t_{1}^{3}}.
\]
\item [{\uline{Type\mbox{$\,$}B}:}] 
\[
\mathbf{g}=-\frac{3}{\Lambda\,t_{1}^{2}}\,\left(dt_{1}^{2}+\epsilon_{0}dt_{2}^{2}\right)+2t_{2}dt_{1}dz_{1}+\frac{\psi}{4t_{1}^{2}}dz_{1}^{2}+\frac{1}{t_{2}^{2}}dz_{1}dz_{2}
\]
with
\[
X=a_{1}\partial_{z_{1}}+\left(a_{2}z_{1}+a_{3}\right)\partial_{z_{2}},
\]
where $\epsilon_{0}\pm1$ and $\Lambda,a_{i}\in\mathbb{R}$ are arbitrary
constants satisfying $\Lambda<0$, moreover $\psi=\psi(t_{1},t_{2})$
is a differentiable function satisfying the differential equation
\[
\psi_{,t_{1}t_{1}}+\epsilon_{0}\psi_{,t_{2}t_{2}}-\frac{2}{t_{1}}\psi_{,t_{1}}=-\epsilon_{0}\frac{4}{3}\Lambda t_{1}^{5}-\frac{12\,a_{2}}{\Lambda\,t_{1}}.
\]
In particular the Ricci soliton is expanding and the orbit metric
$\tilde{\mathbf{g}}=-3\Lambda^{-1}t_{1}^{-2}(dt_{1}^{2}+\epsilon_{0}dt_{2}^{2})$
has constant Gauss curvature $\Lambda/3$.
\end{description}
\end{thm}

\bigskip{}

We notice that the transformation from Type A' to Type A did not preserve
the conformal representation of the orbit metric $\mathbf{\tilde{g}}$,
due to the fact that $\sqrt{R}$ is not harmonic when $c\neq0$. 

Also we notice that the Ricci solitons of previous theorem reduce
to Einstein metrics when $a_{2}=0$. In such a case, by means of the
coordinate transformation $\{t_{1}\mapsto t_{1},\,t_{2}\mapsto t_{2},\,z_{1}\mapsto\,-z_{1},\,z_{2}\mapsto z_{2}-\epsilon\,t_{2}\,/t_{1}^{3}\}$,
the Einstein metric of Type A can be written as
\begin{equation}
\mathbf{g}=\begin{array}{l}
{\displaystyle -\frac{3c_{1}^{2}t_{1}}{\Lambda\left(c_{1}^{2}t_{1}^{3}+1\right)}\,dt_{1}^{2}+\frac{\epsilon_{0}\left(c_{1}^{2}t_{1}^{3}+1\right)}{t_{1}}\,dt_{2}^{2}+\frac{2\,\epsilon}{t_{1}}\,dt_{2}\,dz_{1}+\frac{\epsilon\,t_{1}^{3}\,\psi+\epsilon_{0}}{t_{1}}\,dz_{1}^{2}+2t_{1}^{2}\,dz_{1}\,dz_{2}},\vspace{10pt}\end{array}\label{eq:Kundu_gen1}
\end{equation}
with $\psi=\psi(t_{1},t_{2})$ satisfying 
\begin{equation}
\psi_{,t_{1}t_{1}}-\frac{3\,\epsilon_{0}\,\epsilon\,c_{1}\,t_{1}^{2}}{\Lambda\,(\epsilon\,t_{1}^{3}+c_{1})^{2}}\psi_{,t_{2}t_{2}}+\frac{\left(4\,\epsilon\,t_{1}^{3}+c_{1}\right)}{t_{1}\left(\epsilon\,t_{1}^{3}+c_{1}\right)}\psi_{,t_{1}}=0.\label{eq:EQ-Kundu_gen1}
\end{equation}
Analogously, by means of the coordinate transformation 
\[
t_{1}\mapsto t_{1}\sqrt[3]{\frac{\Lambda}{3}},\qquad t_{2}\mapsto t_{2}\sqrt[6]{-\frac{3}{\Lambda}},\qquad z_{1}\mapsto\frac{z_{1}}{3}\,\sqrt[6]{-\frac{3}{\Lambda}},\qquad z_{2}\mapsto-\frac{\Lambda\,z_{2}}{2}\,\sqrt[6]{-\frac{3}{\Lambda}},
\]
and using $-t_{1}^{6}-6\Lambda\psi$ instead of $\psi$, the Einstein
metric of Type B can be written as 
\begin{equation}
\mathbf{g}=\begin{array}{l}
-{\displaystyle \frac{3}{\Lambda t_{1}^{2}}}\,dt_{1}^{2}+\epsilon_{0}{\displaystyle \frac{1}{t_{1}^{2}}}dt_{2}^{2}+2t_{1}\,dt_{2}\,dz_{1}+{\displaystyle \frac{t_{1}^{6}+\psi}{2t_{1}^{2}}}\,dz_{1}^{2}+{\displaystyle \frac{2}{t_{1}^{2}}}\,dz_{1}\,dz_{2},\end{array}\label{eq:Kundu_gen2}
\end{equation}
with $\psi=\psi(t_{1},t_{2})$ satisfying 
\begin{equation}
\psi_{t_{1}t_{1}}-{\displaystyle \frac{2}{t_{1}}}\psi_{t_{1}}-{\displaystyle \frac{3\epsilon_{0}}{\Lambda}}\psi_{t_{2}t_{2}}=0.\label{eq:EQ-Kundu_gen2}
\end{equation}
The Einstein metrics (\ref{eq:Kundu_gen1}-\ref{eq:EQ-Kundu_gen1})
and (\ref{eq:Kundu_gen2}-\ref{eq:EQ-Kundu_gen2}), with $\epsilon_{0}=1$,
have already been described in the paper \cite{Cat-Mar}, but their
expressions there had some misprints which we amended here. 

\subsection*{Case II-2}

Here is the main result when $f_{1,t_{2}}^{1}=0$ and $h_{12,t_{1}}\neq0$,
$\Lambda=0$.
\begin{thm}
For a suitable choice of adapted coordinates $\{t_{1},t_{2},z_{1},z_{2}\}$,
any Ricci soliton with a $2$-dimensional Abelian Killing algebra
$\mathcal{G}_{2}$, which is orthogonally intransitive, has non-null
Killing leaves and satisfies conditions (i), (ii) and (iii) with $f_{1,t_{2}}^{1}=0$
and $h_{12,t_{1}}\neq0$, $\Lambda=0$, can be written in the form
\[
g=\frac{1}{\sqrt{t_{1}}}\,\left(dt_{1}^{2}+\epsilon_{0}dt_{2}^{2}\right)+2\frac{t_{2}}{t_{1}^{3/2}}\,dt_{1}dz_{1}+\left(\psi\,t_{1}+\frac{4\,\epsilon_{0}}{9\,\sqrt{t_{1}}}\right)\,dz_{1}^{2}+2\,t_{1}dz_{1}dz_{2}
\]
with
\[
X=a_{1}\partial_{z_{1}}+\left(a_{2}z_{1}+a_{3}\right)\partial_{z_{2}},
\]
where $a_{i}\in\mathbb{R}$ are arbitrary constants, $\epsilon_{0}=\pm1$
and in addition $\psi=\psi(t_{1},t_{2})$ is a differentiable function
satisfying the differential equation 
\[
\psi_{,t_{1}t_{1}}+\epsilon_{0}\psi_{,t_{2}t_{2}}+\frac{1}{t_{1}}\psi_{,t_{1}}=\frac{2a_{2}}{\sqrt{t_{1}}}.
\]
In particular the orbit metric $\tilde{\mathbf{g}}=t_{1}^{-1/2}(dt_{1}^{2}+\epsilon_{0}dt_{2}^{2})$
has Gauss curvature $-t_{1}^{-3/2}/4$.
\end{thm}

\bigskip{}

It is easy to check that, when $a_{2}=0$, the Ricci soliton described
by this theorem is equivalent to the Einstein metric found by Kundu
in the paper \cite{Kun}. 

\subsection*{Case II-3}

Here is the main result when $f_{1,t_{2}}^{1}=0$ and $h_{12,t_{1}}=0$.
\begin{thm}
For a suitable choice of adapted coordinates $\{t_{1},t_{2},z_{1},z_{2}\}$,
any Ricci soliton with a $2$-dimensional Abelian Killing algebra
$\mathcal{G}_{2}$, which is orthogonally intransitive, has non-null
Killing leaves and satisfies conditions (i), (ii) and (iii) with $f_{1,t_{2}}^{1}=0$
and $h_{12,t_{1}}=0$, can be written in the form 
\[
g=e^{-2P}\,\left(dt_{1}^{2}+\epsilon_{0}dt_{2}^{2}\right)+2\,c\left(\int e^{-2P}dt_{2}\right)\,dt_{1}dz_{1}+\psi\,dz_{1}^{2}+2\,dz_{1}dz_{2}
\]
with
\[
X=\left(2\Lambda z_{1}+a\right)\partial_{z_{1}}+A\partial_{z_{2}},
\]
where $\epsilon_{0}=\pm1$ and $a,c\in\mathbb{R}$ are arbitrary constants
satisfying $c\neq0$, moreover $A=A(z_{1})$, $P=P(t_{1},t_{2})$
and $\psi=\psi(t_{1},t_{2})$ are differentiable functions satisfying
the differential equations 
\[
\left\{ \begin{array}{l}
\psi_{,t_{1}t_{1}}+\epsilon_{0}\psi_{,t_{2}t_{2}}-2\Lambda e^{-2P}\psi=\epsilon_{0}e^{-2P}\left(c^{2}+2\epsilon_{0}A'\right)\vspace{15pt}\\
P_{,t_{1}t_{1}}+\epsilon_{0}P_{,t_{2}t_{2}}=\Lambda e^{-2P}.
\end{array}\right.
\]
In particular the Ricci soliton can be either expanding ($\Lambda<0$),
steady ($\Lambda=0$) or shrinking ($\Lambda>0$), and the orbit metric
$\tilde{\mathbf{g}}=e^{-2P}\,(dt_{1}^{2}+\epsilon_{0}dt_{2}^{2})$
has constant Gauss curvature $\Lambda$.
\end{thm}

\section*{Acknowledgements}

The author was partially supported by CNPq, grant 310577/2015-2 and
grant 422906/2016-6.

\end{document}